\documentclass[12pt]{article}
\usepackage{amsmath,amssymb,eucal}
\usepackage{hyperref}
\font\tenrm=cmr10

\font\cmssl=cmss10 at 12 pt

\font\bigss=cmssdc10 scaled 2300

\font\cmsslll=cmss10 at 14 pt

\renewcommand{\a}{\alpha}
\renewcommand{\b}{\beta}

\newcommand{\g}{\gamma}

\renewcommand{\i}{\iota}

\renewcommand{\l}{\lambda}
\renewcommand{\o}{\omega}

\renewcommand{\r}{\rho}

\newcommand{\x}{\xi}

\newcommand{\G}{\Gamma}

\newcommand{\bC}{\mathbb{C}}
\newcommand{\bR}{\mathbb{R}}
\newcommand{\bZ}{\mathbb{Z}}

\renewcommand{\gg}{\mathfrak{g}}
\newcommand{\gh}{\mathfrak{h}}
\newcommand{\gk}{\mathfrak{k}}
\newcommand{\gl}{\mathfrak{l}}
\newcommand{\gm}{\mathfrak{m}}

\newcommand{\gs}{\mathfrak{s}}

\newcommand{\gz}{\mathfrak{z}}

\newcommand{\su}{\mathfrak{su}}

\newcommand\GL{\mathrm{GL}}
\newcommand\SL{\mathrm{SL}}
\newcommand\SO{\mathrm{SO}}
\newcommand\SU{\mathrm{SU}}
\newcommand\U{\mathrm{U}}

\newcommand{\p}{\partial}

\renewcommand{\square}{\kern1pt\vbox
               {\hrule height 0.6pt\hbox{\vrule width 0.6pt\hskip 3pt
    \vbox{\vskip 6pt}\hskip 3pt\vrule width 0.6pt}\hrule height0.6pt}
    \kern1pt}

\newcommand{\ra}{\rightarrow}

\DeclareMathOperator\End{End\;}

 \DeclareMathOperator\ad{ad\;}

\newtheorem{Pb}{Problem}

\newtheorem{Th}{Theorem}[section]  
\newtheorem{Prop}[Th]{Proposition}

\newtheorem{Cor}[Th]{Corollary}
\newtheorem{Lem}[Th]{Lemma}
\newtheorem{Def}[Th]{Definition}
\newcommand{\bP}{\begin{Pb}\ \ }
\newcommand{\eP}{\end{Pb}}
\newcommand{\bt}{\begin{Th}\ \ }
\newcommand{\et}{\end{Th}}
\newcommand{\bp}{\begin{Prop}\ \ }
\newcommand{\ep}{\end{Prop}}
\newcommand{\bc}{\begin{Cor}\ \ }
\newcommand{\ec}{\end{Cor}}
\newcommand{\bl}{\begin{Lem}\ \ }
\newcommand{\el}{\end{Lem}}
\newcommand{\bd}{\begin{Def}\ \ }
\newcommand{\ed}{\end{Def}}

\newcommand{\pf}{\noindent{\it Proof:\ \ }}
\newcommand{\qed}{\hfill\square}
\newcommand{\n}{\nabla}

\newcommand{\ot}{\otimes}

\newcommand{\be}{\begin{equation}}
\newcommand{\ee}{\end{equation}}

\newcommand\re[1]{(\ref{#1})}
\newcommand{\arr}{\begin{array}{rlll}}
\newcommand{\ea}{\end{array}}
\newcommand{\bea}{\begin{eqnarray}}
\newcommand{\eea}{\end{eqnarray}}
\newcommand{\bean}{\begin{eqnarray*}}
\newcommand{\eean}{\end{eqnarray*}}

\catcode`@=11
\@addtoreset{equation}{section}
\catcode`@=12

\begin{document}
 \rightline{}
\vskip 1.5 true cm
\begin{center}
{\bigss  Homogeneous locally conformally\\[.5em]
K\"ahler and Sasaki manifolds}\\[.5em]
\vskip 1.0 true cm
{\cmsslll  D.\ V.\ Alekseevsky$^1$, V.\ Cort\'es$^2$, K.\ Hasegawa$^3$ and Y.\ Kamishima$^4$} \\[3pt]
$^1${\tenrm Institute for Information Transmission Problems\\
Bolshoi Karetnyi per., 19,
127994 Moscow, Russia\\
and\\ 
Masaryk University\\
Kotlarska 2, 
61137 Brno, Czech Republic\\
dalekseevsky@iitp.ru}\\[1em] 
$^2${\tenrm   Department of Mathematics\\
and Center for Mathematical Physics\\
University of Hamburg\\
Bundesstra{\ss}e 55,
D-20146 Hamburg, Germany\\
cortes@math.uni-hamburg.de}\\[1em]
$^3${\tenrm Department of Mathematics, 
Faculty of Education\\
Niigata University, 8050 Ikarashi-Nino-cho, Nishi-ku, \\
950-2181 Niigata, Japan\\
hasegawa@ed.niigata-u.ac.jp}\\[1em]
$^4${\tenrm Department of Mathematics, Josai University\\
Keyaki-dai 1-1, Sakado, 
350-0295 Saitama, Japan\\
kami@josai.ac.jp}\\[1em]
\end{center}
\vskip 1.0 true cm
\baselineskip=18pt
\begin{abstract}
\noindent
We prove various classification results for homogeneous 
locally conformally symplectic manifolds. 
In particular, we show that a homogeneous locally conformally 
K\"ahler manifold of  a reductive group is of Vaisman type 
if the normalizer of the isotropy group is compact. 
We also show that such a result does not hold 
in the case of non-compact normalizer and  determine 
all left-invariant locally conformally K\"ahler structures 
on reductive Lie groups. 
\end{abstract}
\break

\tableofcontents

\section*{Introduction}
Recall that the notion of a locally conformally K\"ahler manifold $(M,\o , J)$ is 
a generalization of the geometric structure encountered on the Hopf manifolds \cite{V}, see Definition \ref{lcKDef}. 
The study of locally conformally K\"ahler manifolds goes beyond the framework of K\"ahler and symplectic geometry
while still remaining within that of complex and Riemannian geometry. Ignoring the 
complex structure $J$ one arrives at the more general notion of a locally conformally symplectic manifold $(M,\o)$. 
Such manifolds were first considered in \cite{L}. The fundamental $2$-form $\o$ satisfies the 
equation 
\[ d\o = \l \wedge \o \]
for some closed $1$-form $\l$, see Definition \ref{lcsDef}. 
The relation between 
locally conformally K\"ahler manifolds and locally conformally symplectic manifolds
is analogous to the one existing between K\"ahler manifolds and symplectic manifolds. 

This work started in September 2010 during a meeting in Japan with discussions about  the work of 
Hasegawa and Kamishima on compact homogeneous 
locally conformally K\"ahler manifolds. And conversely, some of the results of this collaboration have influenced \cite{HK1}
and \cite{HK2} where the present paper is referenced. This applies in particular to the proof of 
Theorem \ref{lcKThm} that a homogeneous locally conformally K\"ahler manifold of  a reductive group is
of Vaisman type if the normalizer of the isotropy group is compact. In the special case of compact groups,
this theorem has been proved in \cite{HK2} and \cite{GMO} (c.f.\ \cite{MO} for a proof under additional assumptions).  

Now we describe the structure of this article and mention some of its main results. 
In the first section we describe some general constructions relating 
sympletic manifolds, contact manifolds, symplectic cones and locally conformally
symplectic manifolds. In the second section we prove more specific results 
relating K\"ahler manifolds, Sasaki manifolds, K\"ahler cones and 
locally conformally K\"ahler manifolds. The main new object is 
an integrable complex structure compatible with the geometric structures 
considered in the first section. We believe that the systematic presentation 
in the first two sections of the paper 
is useful although part of the material is certainly known to experts in the field. 
In any case,  it is a basis for our investigation of  
homogeneous locally symplectic and locally conformally K\"ahler manifolds 
in the third and fourth sections respectively. Under rather general assumptions,
we first prove that the dimension of the center of a Lie group of automorphisms of a 
locally conformally symplectic manifold
is at most $2$. The main result of the third section
is then a classification of all homogeneous locally symplectic manifolds $(M=G/H,\o)$ with
trivial twisted cohomology class 
$[\o ] \in H^2_\l (\gg , \gh )$ (see Theorem \ref{mainlcsThm}). 
These assumptions are satisfied if $\gg$ is reductive (see Proposition \ref{cohProp}). 

In the last and main section we focus on homogeneous locally
conformally K\"ahler manifolds of reductive groups. As a warm up, we begin by 
classifying left-invariant locally conformally K\"ahler structures on four-dimensional 
reductive Lie groups. We find that not all of them are of Vaisman type.
In Theorem \ref{linvThm} we give the classification of left-invariant locally conformally K\"ahler structures
on arbitrary reductive Lie groups. The case of general homogeneous spaces $G/H$ 
of reductive groups $G$ is related to the case of trivial stabilizer $H$ by considering
the induced locally conformally K\"ahler structure on the Lie group $N_G(H)/H$. 
Assuming the latter group to be compact, we prove that the initial locally conformally K\"ahler structure 
on $G/H$ is necessarily of Vaisman type (see Theorem \ref{lcKThm}). 
  
 \noindent
{\bf Acknowledgments}\\ 
D.V.A.\ has been supported by the project CZ.1.07/2.3.00/20.0003 of the Operational
Programme Education for Com\-petitiveness of the Ministry of Education, Youth and Sports 
of the Czech Republic. V.C. has been supported by the RTG 1670 
``Mathematics inspired by String Theory'', funded by the Deutsche Forschungsgemeinschaft (DFG). 
K.H.\ and Y.K.\ have been supported by JSPS Grant-in-Aid for Scientific Research. 

D.V.A.\ thanks the University of Hamburg, Niigata University
and Tokyo Metropolitan University for hospitality and support. 
V.C.\ thanks Masaryk University (Brno), Niigata University 
and Tokyo Metropolitan University for hospitality and support. 
K.H.\ and \mbox{Y.K.}\ thank the University of Hamburg and Masaryk University for hospitality and support. 

We would like to thank very much the referee for useful comments and remarks.


\section{
Symplectic manifolds,
contact manifolds and symplectic cones}

\subsection{Contactization}
\bd A symplectic manifold $(M,\o )$ is called {\cmssl A-quantizable} if there
exists a principal bundle $\pi : P\ra M$ with
one-dimensional structure group $A= S^1$ or $\bR$
and connection $\theta$ such that $d\theta= \pi^*\o$.
\ed
The closed $2$-form $\o$ gives rise to a \v{C}ech cohomology class
$[c]\in \check{H}^2(M,\bR)$, which can be defined as follows.
Let $(U_\a)$ be a covering of $M$ by contractible open
sets such that the intersections $U_{\a\b}:= U_{\a}\cap U_{\b}$
and $U_{\a\b\g} :=U_{\a}\cap U_{\b}\cap U_{\g}$
are also contractible. By the Poincar\'e Lemma, on each $U_\a$ we can choose a
$1$-form $\theta_\a$ such that $d\theta_\a = \o|_{U_\a}$.
Similarly, the $1$-form
\[\theta_{\a\b} := \theta_\a|_{U_{\a\b}}-\theta_\b|_{U_{\a\b}}\]
is closed and, hence, $\theta_{\a\b} = df_{\a\b}$ for some function $f_{\a\b}
=-f_{\b\a}\in C^{\infty}(U_{\a\b})$. Finally, the function
\[ c_{\a\b\g} := f_{\a\b}|_{U_{\a\b\g}}+  f_{\b\g}|_{U_{\a\b\g}}+ f_{\g \a}|_{U_{\a\b\g}}\]
is closed  and hence constant. By construction, $c = (c_{\a\b\g})$ is a
\v{C}ech 2-cocycle with values in the constant sheaf $\bR$. One can check that
the corresponding class $[c]\in \check{H}^2(M,\bR)$
depends only on the de Rham cohomology class
$[\o]\in H^2(M,\bR)$. We will call $[c]$ the {\cmssl characteristic class}
of the symplectic manifold $(M,\o)$.
Recall that a class $[c]\in \check{H}^2(M,\bR)$ is called {\cmssl integral}
if it can be represented by an integral cocycle, that is a
cocycle $c=(c_{\a\b\g})$ such that $c_{\a\b\g}\in \bZ$.
\bp A symplectic manifold $(M,\o )$ is $S^1$-quantizable if and only if
its characteristic class $[c]\in \check{H}^2(M,\bR)$ is integral.
It is  $\bR$-quantizable if and only if $[c]=0$. In particular,
any exact symplectic manifold is quantizable.
\ep

\bd Any such pair $(P,\theta )$ will be called a {\cmssl contactization}
(or, more precisely, {\cmssl A-contactization}, where $A=S^1$ or $\bR$)
of the symplectic manifold $(M,\o)$. By a {\cmssl contact manifold}
we will understand a manifold $P$ of dimension $2n+1$
together with a globally defined contact
form $\theta$, that is $d\theta^n \wedge \theta\neq 0$.
A contact manifold $(P,\theta )$ will be
called {\cmssl regular} if its Reeb vector field $Z$ generates a free and
proper action of $A=S^1$ or $\bR$.
\ed

\bp Any contactization $(P,\theta )$ of an $A$-quantizable
symplectic manifold $(M,\o )$ is a regular contact manifold
with global contact form $\theta$. The group $\mathrm{Aut} (P,\theta )$ contains
the 1-dimensional central subgroup $A$, which is the kernel of
the natural homomorphism  $\mathrm{Aut} (P,\theta ) \ra \mathrm{Aut}(M,\o )$.
\ep

\pf $\theta$ is indeed a contact form, since
$d\theta=\pi^*\o$ is non-degenerate on the
horizontal distribution $\ker \theta$. The Reeb vector field $Z$ is the
generator of the principal action, which is free and proper. \qed

\bp
There is a bijection between $A$-quantizable symplectic manifolds $(M,\o)$ with $H^1 (M, \bR) =0$ up to isomorphism
and regular contact manifolds $(P,\theta )$ with Reeb action of $A=S^1$ or $\bR$ up to isomorphism.
\ep

\subsection{Symplectic cone over a contact manifold}
Let $(P,\theta )$ be a contact manifold. We denote by
$N = C(P) = \bR^{>0}\times P$ the cone over $P$ with the
radial coordinate $r$.
\bp For any contact manifold $(P,\theta )$,
\[ \o_N := rdr\wedge \theta + \frac{r^2}{2}d\theta = d(\frac{r^2}{2}\theta )\]
is a symplectic form on the cone $N=C(P)$.
\ep
\bd \label{symplconeDef} 
The pair $(N,\o_N)$ is called {\cmssl the symplectic cone}
over the contact manifold $(P,\theta )$.
\ed
Now we give an intrinsic characterization of symplectic cones
in the category of symplectic manifolds.
\bd A {\cmssl conical symplectic manifold} $(M,\o ,\xi, Z)$
is a symplectic manifold $(M,\o )$ endowed with two commuting vector fields
$\xi$ and $Z$ such that
\[ \o (\xi , Z) >0,
\quad \mathcal{L}_\xi \o = 2 \o,\quad \mathcal{L}_Z \o = 0.\]
A {\cmssl global
conical symplectic manifold} is a conical symplectic manifold
$(M,\o ,\xi, Z)$ such that $\xi$ is complete.
\ed
\bt \begin{enumerate}
\item[(i)] The symplectic cone over any contact manifold is a global conical
symplectic manifold.
\item[(ii)] Conversely, any global conical symplectic manifold is a
symplectic cone over a contact manifold.
\item[(iii)] Any conical symplectic manifold is locally isomorphic
to a symplectic cone over a contact manifold.
\end{enumerate}
\et

\pf (i) Let $(N=C(P),\o_N)$ be a symplectic cone over a contact manifold
$(P,\theta )$. The Reeb vector field of $P$ can be considered
as a vector field $Z$ on $N$, which together with $\xi = r\p_r$
defines a global conical structure. To prove (ii-iii)
we need the following lemma.

\bl Let $(M,\o ,\xi, Z)$ be a conical symplectic manifold. Let
$f$ be a positive smooth function defined in some
open neighborhood $U$ such that $df=-\i_Z\o$, i.e.\
$f$ is the Hamiltonian of $-Z$. Then in $U$ the symplectic form
$\o$ can be written as
\[ \o= df\wedge \theta + fd\theta = rdr\wedge \theta + \frac{r^2}{2}d\theta,\]
where
\[ \theta = \frac{1}{2f}\eta,\quad \eta = \i_\xi\o,\quad r= \sqrt{2f}.\]
\el

\noindent
{\bf Remark:} The function $f$ is unique up to addition of
a constant $c$ such that $f+c>0$. We can choose,
for example, $f=\frac{1}{2}\o (\xi , Z)$, which is characterized
by the condition $\mathcal{L}_\xi f=2f$.

\pf The symplectic form is exact:
\[ 2\o = d\eta ,\quad \eta:= \i_\xi\o.\]
We define
\[ \theta := \frac{1}{2f}\eta.\]
Then we calculate
\[ df\wedge \theta + fd\theta = \frac{df}{2f}\wedge \eta + fd(\frac{1}{2f})
\wedge \eta + \o = \o.\]
Now it suffices to rewrite
\[ f=\frac{r^2}{2}\]
to obtain $\o = rdr\wedge \theta + \frac{r^2}{2}d\theta$. \qed

The lemma proves part (iii) of the theorem. To prove (ii)
we remark that using the flow of the complete vector field
$\xi$ on a global conical symplectic manifold $(N,\o ,\xi,Z)$ we get a global diffeomorphism $N\cong I\times P$,
where $P$ is some level set of $f=\frac{1}{2}\o (\xi , Z)$
and $I = (a,b)$, where
$0\ge a=\inf f$, $b=\sup f$. We have to show that $a=0$ and $b=\infty$.
Let $\g : \bR \ra N$ be an integral curve of $\xi$. Then
$\mathcal{L}_\xi f=2f$ implies the differential equation
$h'=2h$, where $h=f\circ \g$. Therefore,
$h(t)=c e^{2t}$ for some positive constant $c$, since $f>0$.
This shows that $I=\bR^{>0}$ and that $N$ is a
symplectic cone $N=C(P)$, where $P=\{ r=1 \} = \{f=1/2\}$. \qed

\subsection{Symplectic cones and locally conformally symplectic
manifolds}
\bd \label{lcsDef} 
A {\cmssl locally conformally symplectic manifold  (lcs manifold) } $(M,\o)$ is a
smooth manifold endowed with a non-degenerate
$2$-form such that $d\o = \l \wedge \o$ for some
closed $1$-form $\l$ called {\cmssl Lee form}. 
An lcs manifold is called {\cmssl proper} if $d\o \neq 0$. 
The vector field $Z:=
\frac{1}{2}\o^{-1}\l$
is called the {\cmssl Reeb field}.
\ed
\noindent
{\bf Remark:} Since $\o$ is non-degenerate,
the equation $d\o = \l \wedge \o$ implies
$d\l=0$ provided that $\dim M >4$.
\bp \label{autoProp}
The vector field $Z$ is an infinitesimal
automorphism of   $(M,\o)$.
\ep
\pf
\[ \mathcal{L}_Z\o = d\i_Z\o+\i_Zd\o = \frac{1}{2}d\l + \i_Z(\l \wedge\o) =0,\]
since $\l (Z) = 2\o (Z,Z)=0$ and $\l\wedge \l =0$.
\qed

Let $(N,\o_N)$ be a symplectic cone over a contact manifold $(P,\theta)$. We define
\[ \o_{lcs} := \frac{1}{r^2}\o_N = dt\wedge \theta + \frac{1}{2}d\theta,\quad
t=\ln r. \]
\bp \label{lcsProp}
For any non-trivial discrete subgroup $\G \subset \bR^{>0}$
the manifold
$(N/\G=S^1\times P,\o_{lcs})$ is lcs.
\ep

\section{
K\"ahler manifolds, Sasaki manifolds and K\"ahler cones}
\subsection{Contactizations of K\"ahler manifolds}
\bd \label{SasakiDef} A {\cmssl Sasaki manifold} $(S,g,Z)$ is a Riemannian
manifold $(S,g)$ endowed with a unit Killing vector
field $Z$, such that $J:=\n Z|_{\mathcal{H}}$  defines an integrable
CR structure on the distribution $\mathcal{H} := Z^\perp\subset TS$.
\ed
Let $(S,g,Z)$ be a Sasaki manifold. Then we define the
$1$-form
\[ \theta := g(Z,\cdot ).\]

\bp For any Sasaki manifold $(S,g,Z)$ the $1$-form $\theta$
is a contact form with the Reeb vector field $Z$
and the CR structure is strictly pseudo-convex.
\ep

\pf It follows from Definition \ref{SasakiDef} that $d\theta =
g(J\cdot ,\cdot )$  on $Z^\perp=\ker \theta$ is non-degenerate.
Hence, $\theta$ is a contact form with positive definite Levi form.
Furthermore,  $\theta (Z)=1$ and
\[ 0=\mathcal{L}_Z\theta = \i_Zd\theta ,\]
which shows that $Z$ is the Reeb vector field.
 \qed

The following theorem establishes a one-to-one correspondence between
quantizable K\"ahler manifolds and
regular Sasaki manifolds.
\bt  Let $A=S^1$ or $\bR$.
\begin{enumerate}
\item[(i)]
The contactization of an $A$-quantizable K\"ahler manifold
$(M,\o , J)$ is a regular Sasaki manifold $(S, \theta , g_S, Z)$,
where $(S,\theta )$, $\pi : S \ra M=S/A$, is the contactization of
$(M,\o )$ with the
fundamental vector field $Z$ of the $A$-action and
\[ g_S = \theta^2 + \frac{1}{2}\pi^*g_M,\quad g_M = \o (\cdot ,J\cdot ).\]
\item[(ii)]
Conversely, any regular Sasaki manifold
with Reeb action
of $A$ is the contactization of an $A$-quantizable K\"ahler manifold.
\end{enumerate}
\et

\subsection{Cones over Sasaki manifolds and K\"ahler cones}

\bd A {\cmssl conical Riemannian manifold} $(M,g,\xi )$ is a
Riemannian manifold $(M,g)$ endowed with a nowhere vanishing
(homothetic) vector field $\xi$ such that $\n \xi = \mathrm{Id}$.
If $\xi$ is complete it is called a  {\cmssl global conical
Riemannian manifold}.
\ed

\bp \begin{enumerate}
\item[(i)] The metric cone over any Riemannian manifold is a global conical
Riemannian manifold.
\item[(ii)] Conversely, any global conical Riemannian manifold is a
metric cone.
\item[(iii)] Any conical Riemannian manifold is locally isometric
to a metric cone.
\end{enumerate}
\ep

\bd A {\cmssl K\"ahler cone} $(N,g_N,J)$ is a
metric cone $(N=C(M),g_N=dr^2+r^2g_M)$ over a Riemannian
manifold $(M,g_M)$ endowed with a skew-symmetric parallel complex
structure $J$.
\ed

\bp Any conical K\"ahler manifold is locally a K\"ahler cone
and any global conical K\"ahler manifold is a K\"ahler cone.
\ep

\bt \label{KaehlerconeThm} \begin{enumerate}
\item[(i)] The metric cone $(N=C(S),g_N)$ over a Sasaki manifold $(S,g_S,Z)$
equipped with
the complex structure $J_N$ defined by
\[ J_N|_{\mathcal{H}} := J = \n Z|_{\mathcal{H}},\quad J_N\xi := Z,\]
is a K\"ahler cone.
\item[(ii)] Conversely, any K\"ahler cone is the cone
over a Sasaki manifold and any conical K\"ahler manifold
is locally isomorphic to a K\"ahler cone over a Sasaki manifold.
\end{enumerate}
\et

Now we give a characterisation of Sasaki manifolds in the class
of strictly pseudo-convex CR manifolds. In the same way one can
characterize pseudo-Riemannian Sasaki manifolds in the class of
Levi non-degenerate CR-manifolds. 

Let $(P,\theta ,J )$ be a strictly pseudo-convex 
integrable CR-structure with globally defined contact form $\theta$, 
which defines the (contact) CR-distribution 
$\mathcal{H}=\ker \theta$. We denote by $Z$ the Reeb vector field of 
$\theta$, such that $\theta (Z) =1$ and $d\theta (Z,\cdot )=0$
and extend $J$ defined on $\mathcal{H}$ to an endomorphism field
on $TP= \bR Z \oplus \mathcal{H}$ by $JZ=0$. Then we 
define a natural Riemannian metric 
$g_P$ on $P$ by 
\[ g_P := \theta^2 + \frac{1}{2}d\theta (\cdot , J\cdot ).\] 
The vector field $Z$ preserves $\theta$ but does
not preserve $J$ and $g_P$ in general. 

\bt Let $(P,\theta ,J)$ be a strictly pseudo-convex 
integrable CR-structure with globally defined contact form $\theta$. 
Then the symplectic structure  $\o_N$ of the symplectic cone 
$(N,\o_N)$ over the contact manifold $(P,\theta )$ (see 
Definition  \ref{symplconeDef}) together with the cone metric
$g_N=dr^2+r^2g_P$ defines on $N=C(P)= \bR^{>0}\times P$ an
almost K\"ahler structure. It is K\"ahler if and only if
the Reeb vector field is {\cmssl holomorphic}, that is 
an infinitesimal CR-automorphism: $\mathcal{L}_ZJ=0$.   
\et 

\pf We have to check that the skew-symmetric endomorphism
$J_N = g_N^{-1}\circ \o_N$ is an almost complex structure. 
Recall that 
\begin{eqnarray*} \o_N &=&  rdr\wedge \theta +  \frac{r^2}{2}d\theta ,\\ 
g_N&=& dr^2 + r^2\theta^2 + \frac{r^2}{2}d\theta (\cdot , J\cdot ).
\end{eqnarray*}
{}From these formulas we see that the decomposition $\mathcal{H} \oplus 
\mathrm{span} \{ \partial_r,Z\}$ is orthogonal with respect to
$\o_N$ and $g_N$. Hence, $J_N$ preserves this 
decomposition and $J_N|_{\mathcal{H}}=J$. We check that $J_NZ = -\xi 
:=-r\partial_r$ and $J_N\xi = Z$:
\begin{eqnarray*}&& \o_N(Z,\cdot ) = -rdr = -g_N(\xi ,\cdot ), \\
&&\o_N (\xi ,\cdot ) = r^2\theta = g_N( Z,\cdot ).  
\end{eqnarray*}
Now we investigate the integrability of $J_N$, that is the 
involutivity of $T^{0,1}N\subset T^{\bC}N$. The involutivity
of $\mathcal{H}^{0,1}$ follows from the integrability
of the CR-structure $J=J_N|_{\mathcal{H}}$. The involutivity
of $(\mathcal{H}^\perp)^{0,1}=\bC (Z+i\xi )$ is automatic
for dimensional reasons. Finally the bracket of $Z+iJ_NZ=Z-i\xi$ with 
$X+iJ_NX=X+iJX$, $X\in \G (P, \mathcal{H})\subset \G (N,\mathcal{H})$, 
is computed as follows:
\[ [Z+i\xi , X+iJX] = [Z,X+iJX] = [Z,X] +i [Z,JX],\]
which is of type $(0,1)$ if and only if 
$[Z,JX] = J[Z,X]$ for all $X$, that is if and only if
$\mathcal{L}_ZJ=0$. 
\qed 

As a corollary,  cf.\ Theorem \ref{KaehlerconeThm}, 
we obtain the following (connection-free) characterization
of Sasaki manifolds in terms of CR-structures.

\bc A Sasaki manifold $(P,g,Z)$ is the same as a 
strictly pseudo-convex CR-manifold $(P,\theta,J)$  
with globally defined contact form  $\theta$ 
such that the corresponding Reeb vector field
$Z$ is holomorphic. The metric $g=g_P$ is the natural Riemannian
metric on $P$ defined by the data $(\theta ,J)$. 
\ec 

\bt \label{SasakiThm} Let $(S_i,g_i,Z_i)$, $i=1,2$, be two Sasaki manifolds.
Then the manifold $N=S_1\times S_2$ has a two-parameter family
of integrable complex structures $J=J_{a,b}$ defined by
\[ J|_{\mathcal{H}_i} = J_i,\quad JZ_1 = aZ_1+ bZ_2,\quad
JZ_2=cZ_1-aZ_2,\]
where $a\in \bR$, $b\neq 0$, $c=-\frac{1+a^2}{b}$ and $(\mathcal{H}_i,J_i)$
is the CR structure of $S_i$. The complex structures
$J_{can}:=J_{0,1}$ and $-J_{can}:=J_{0,-1}$ are
the only structures in the family $J_{a,b}$ for which the
product metric is Hermitian.
\et

\pf This follows from the Newlander-Nirenberg theorem by a direct calculation. 
\qed 

As a special case we obtain the famous complex structures on products of spheres,
constructed by Calabi and Eckmann.
\bc The product of two odd-dimensional spheres has a two-parameter family
$J_{a,b}$ of integrable complex structures. The
product metric is Hermitian with respect to the complex structure
$J_{can}$.
\ec

\subsection{K\"ahler cones and locally conformally K\"ahler
manifolds}
\bd \label{lcKDef} A {\cmssl locally conformally K\"ahler manifold} (lcK manifold) $(M,\o ,J)$
is a locally conformally symplectic manifold $(M,\o)$ endowed
with a skew-symmetric integrable complex structure $J$ such that
the metric
\[ g = \o (\cdot ,J\cdot )\]
is positive definite. The Riemannian metric $g$ is then called a
{\cmssl locally conformally K\"ahler metric} (lcK metric).
The $1$-form $\theta := \frac{1}{2}J^*\l$ is called the 
{\cmssl Reeb form}. 
The (locally gradient) vector field $\xi = -\frac{1}{2}g^{-1}\l$
is called the {\cmssl Lee field}.
An lcK manifold $(M,\o ,J)$ is called {\cmssl Vaisman manifold}
if $\xi$ is a parallel unit vector field. 
\ed

\noindent
Remark that if $\xi$ is parallel then $\l (\xi )$ is constant. By rescaling
$\o$ we can always normalize
$\l (\xi ) =2\o (Z,\xi ) = 2g(JZ ,\xi ) = -2g(\xi ,\xi)= -2$, such that
$|\xi |=1$. Note that, as a consequence of the above definition,  the Lee and the Reeb field are 
related by 
\[ Z=J\xi .\] 

Similarly one defines the notion of a locally conformally pseudo-K\"ahler manifold
and that of a pseudo-Riemannian Vaisman manifold 
by allowing the metric to be indefinite. 

Vaisman manifolds were first studied by Vaisman, who called them
\emph{generalized Hopf manifolds}. In \cite{V} he proved the following
theorem, which relates them to Sasaki manifolds. 
For convenience of the reader we reprove it within the logic 
of our exposition.  
\bt \label{VaismanThm} Let $(M,\omega ,J)$ be a complete Vaisman manifold. Then 
\begin{enumerate}
\item[(i)] the Lee field $\xi$ and the Reeb field $Z=J\xi$ are infinitesimal
automorphisms of the lcK structure $(\omega ,J)$ and 
\item[(ii)] the universal 
cover of $M$ is a Riemannian product of a line
and a simply connected Sasaki manifold $S$.
\end{enumerate}
\et

\pf
The de Rham theorem implies that the universal cover of a
complete Vaisman manifold is a Riemannian product $M=\bR \times S$
of a line and a simply connected manifold $S$, where $S$ is a leaf
of the integrable distribution $\ker \l = \xi^\perp$.
We already know that $\xi$ is
a Killing vector field, since it is parallel. We also know that 
$Z$ preserves $\o$ by Proposition \ref{autoProp}. Therefore, in order
to prove (i), we only have to show that $\xi$ and $Z$ are holomorphic, that is
preserve the complex structure $J$. We recall that a (real) vector field 
$X$ is holomorphic if and only if $JX$ is
holomorphic. Moreover, under this assumption, $X$ and $JX$ commute.
Since $Z=J\xi$, it suffices to check that $\xi$ is holomorphic. 
Now any lcK manifold $(M,\o ,J)$ admits a canonical torsion-free
complex connection $\tilde{\n}$, which coincides with the Levi-Civita
connection of the locally defined K\"ahler metric $\tilde{g}=e^{-f}g$,
where $f$ is a locally defined function such that $df =\l$. 
Indeed, since $f$ is unique up to an additive constant, 
the metric $\tilde{g}$ is unique up to a constant factor and its 
Levi-Civita connection is a well defined connection on $M$.  
With our conventions, the explicit expression for $\tilde{\n}$ is
\be \label{nablatildeequ}
\tilde{\n}_XY = \n_XY -\frac{1}{2}\l (X)Y -\frac{1}{2}\l (Y)X
-g(X,Y)\xi .\ee   
To prove this formula, it is enough to check that the torsion-free connection  
on the right hand side preserves the metric $\tilde{g}$. This is 
a straighforward calculation. 
Using $\n \xi =0$ and \re{nablatildeequ}, we obtain  
$\mathcal{L}_\xi J=  \n_\xi J = \tilde{\n}_\xi J=0$, as in \cite{V}. 

 It follows from (i) that
$\mathcal{L}_\xi \theta =0$. This means that $\theta$ can be considered
as a $1$-form on $S$. 
\bl \label{LxioLemma} Let $(M,\o ,J)$ be an lcK manifold. Then
\[ \mathcal{L}_\xi\o = \l (\xi ) \o -\l \wedge \theta +d\theta.\]
\el

\pf
We calculate
\[  \mathcal{L}_\xi\o = d\theta + \i_\xi (\l\wedge \o)=
d\theta + \l (\xi)\o -\l \wedge \theta .\]
\qed

\noindent
Under the assumptions of the theorem we have $\l (\xi )=-2$, $\theta (Z) =1$
and $\mathcal{L}_\xi \o =0$ such that
\[ \o = -\frac{1}{2}\l \wedge \theta + \frac{1}{2}d\theta.\]
This implies that $d\theta|_S = 2\o|_S$ has 1-dimensional kernel $\bR Z$
transversal to $\mathcal{H}=\ker \theta = Z^\perp$. We have shown that
$\theta$ is a contact form on $S$ with Reeb vector field $Z$. In order
to prove that $S$ is Sasakian, we choose a local function $t$ such that
$\l = -2dt$. Then we can rewrite $\o$ and $g$ in the form
\begin{eqnarray*} \o &=& dt\wedge \theta +\frac{1}{2}d\theta\\
g&=& dt^2 + \theta^2 + \frac{1}{2}\bar{g},
\end{eqnarray*}
where
\be \label{LeviEqu} \bar{g} = d\theta (\cdot ,J\cdot )\ee
is the Levi form.
One can easily check that the metric $g_K=e^{2t}g$ is
a K\"ahler metric with K\"ahler form $\o_K=e^{2t}\o=d(\frac12 e^{2t}\theta )$.
The substitution $r= e^t$ yields
\[ g_K = dr^2 + r^2g_S,\quad g_S = \theta^2 + \frac{1}{2}\bar{g},\quad
\xi = \p_t = r\p_r.\]
This is locally a K\"ahler cone and, hence, its covariant derivative
$\n^K$ yields
\[ \n^K \xi = \mathrm{Id},\quad \n^K Z = \n^K (J\xi) = J.\]
Notice that $g_K|_S$ and $g_S$ are homothetic and, hence,
the Levi Civita connection $\n^S$ of $(S,g_S)$ coincides with
the connection induced by $\n^K$ on the totally umbilic submanifold
$S\subset (M,g_K)$. {}From the Gau{\ss} equation we get
\[  \n^S_XZ= JX\quad\mbox{for all}\quad X\in TS\cap Z^\perp,\quad
\n^S_ZZ =0.\]
This proves that $(S,g_S,Z)$ is a Sasaki manifold.
\qed

\noindent
{\bf Remark:} The isometry group of a compact Vaisman manifold
does not necessarily preserve the complex structure. It suffices
to consider $S^1\times S^{2n+1}$ endowed with the product
metric and the complex structure $J_{can}$ of Theorem \ref{SasakiThm}.
This is an example of an lcK manifold as shown in the next proposition.

Let $(N,\o_N,J_N)$ be a K\"ahler cone over a
Sasaki manifold $(S,g_S,Z)$. Recall that $\o_{lcs}
= dt\wedge \theta + \frac{1}{2}d\theta$
is a conformally symplectic structure on $N$, where
$\theta = g(Z,\cdot )$ is the contact form and $t=\ln r$.
\bp \label{VaismanProp}
For any non-trivial discrete subgroup $\G \subset \bR^{>0}$
the complex structure $J_N$ on the K\"ahler cone $N$ induces
a complex structure $J$ on $N/\G=S^1\times S$ such that
$(N/\G,\o_{lcs},J)$ is a Vaisman manifold.  The group
$S^1=\bR^{>0}/\G$ acts freely, holomorphically and
isometrically (with respect to the
lcK metric) on the lcK manifold $N/\G$ and $Z$ is an $S^1$-invariant
holomorphic Killing vector field on $N/\G$.
\ep

\pf By Proposition \ref{lcsProp}, $(N/\G,\o_{lcs})$ is locally
conformally symplectic. Therefore to prove that it is lcK
it suffices to show that $J_N$ is
invariant under the group $\bR^{>0}$ and, hence, induces a complex structure
$J$ on $N/\G$. This follows from the
equations $\mathcal{L}_\xi\o_N=2\o_N$, $\mathcal{L}_\xi g_N=2g_N$,
since $J_N=g_N^{-1}\o_N$. The group $\bR^{>0}$ acts isometrically
on $N$ with respect to the Riemannian metric
\be  \label{lcKEqu} \o_{lcs}(\cdot , J_N\cdot ) = dt^2+g_S,\ee
which induces the lcK metric $g_{lcK}$ on $M$. In fact $\xi = \partial_t$
is an obvious Killing vector field for the metric \re{lcKEqu}.
This shows that $S^1$ acts isometrically on $(N/\G,g_{lcK})$.
Obviously $\xi =\p_t$ is a parallel unit field and preserves the $2$-form
$\o_{lcs}=dt\wedge \theta + \frac{1}{2}d\theta$. In particular,
$(N/\G,\o_{lcs},J)$ is a Vaisman manifold.
\qed

\noindent
The above complex structure on $N/\G=S^1\times S$ coincides with
the complex structure $J_{can}$ of Theorem \ref{SasakiThm}.
The next theorem shows that the Vaisman manifolds of Proposition
\ref{VaismanProp} admit a canonical two-parameter family of Vaisman deformations.

\bt Let $(N=\bR^{>0}\times S,\o_N,J_N)$ be a K\"ahler cone over a
Sasaki manifold $(S,g_S,Z)$ endowed with the locally
conformally symplectic structure $\o_{lcs}
= dt\wedge \theta + \frac{1}{2}d\theta$. Then
 $(\o_{lcs},J_{a,b})$, where $J_{a,b}$ is defined in Theorem \ref{SasakiThm},
is a Vaisman lcK structure on $N/\G=S^1\times S$
if and only if $b>0$.  
The Reeb vector field $Z$ and the Lee vector field $\xi_{a,b}=-J_{a,b}^*Z$
are holomorphic Killing vector fields
for all of these structures.
\et

\pf $J_{a,b}$ is skew-symmetric with respect to $\o_{lcs}$, since
\[ g_{a,b} := -\o_{lcs}(J_{a,b}\cdot ,\cdot ) =
bdt^2 -c\theta^2 -2adt\theta +\frac{1}{2}\bar{g} \]
is symmetric. (Recall that $\bar{g}$ stands for the Levi form
of $S$, see  \re{LeviEqu}). The metric $g_{a,b}$ is positive definite
if and only if $b>0$. Since $J_{a,b}$ is integrable, by Theorem \ref{SasakiThm},
we see that $(S^1\times S,\o_{lcs},J_{a,b})$ is lcK if $b>0$. 
The vector fields $\xi_{can}=\xi_{0,1}=\p_t$ and $Z$ preserve the 
$1$-forms $dt$ and $\theta$ and, hence, the
metrics $g_{a,b}$. Since the Reeb field always preserves $\o$, this
implies that both vector fields are holomorphic for all $J_{a,b}$.
As a consequence, any linear combination of $\p_t$ and $Z$, such as
$\xi_{a,b}$, is also a holomorphic Killing vector field for any of the
complex structures in the two-parameter family. 
It remains to check
that the lcK structure $(\o_{lcs},J_{a,b})$ is Vaisman. 
The Lee field $\xi_{a,b}=-\frac{1}{2}g_{a,b}^{-1}\l$
is given by
\[ \xi_{a,b} = -c\p_t +aZ.\]
A direct calculation using the Koszul formula for $g=g_{a,b}$
shows that for all
$X, Y \in \mathcal{H}=\ker \theta \cap \ker \l \subset TN$ we have
\begin{eqnarray*} 2g(\n_XY,\p_t) &=& g([X,Y],\p_t) = -a \theta ([X,Y])\\
2g(\n_XY,Z) &=& -Zg(X,Y) + g([X,Y],Z) -g(X,[Y,Z])-g(Y,[X,Z]) =
-c\theta ([X,Y]),
\end{eqnarray*}
since $\mathcal{L}_Zg=0$. As consequence, we obtain
\[ g(\n_X \xi_{a,b},Y) = -g(\n_XY,\xi_{a,b}) = \frac{1}{2} (ac -ca)
\theta ([X,Y])=0,\]
for all $X, Y  \in \mathcal{H}$.  Using the fact that
$\xi_{a,b}$ is a holomorphic Killing vector field, proven above, 
we see that to prove $\n \xi_{a,b}=0$ it is enough to check that 
$\n_{\xi_{a,b}}\xi_{a,b}\perp \mathcal{H}$. Let $X\in \Gamma (\mathcal{H})$ be a local section, which 
commutes with $\xi_{a,b}$. Then the Koszul formula yields   
\[ 2g(\n_{\xi_{a,b}} \xi_{a,b} , X) = -Xg(\xi_{a,b},\xi_{a,b})=0.\] \qed

\bc The Vaisman manifold $(S^1\times S^{2n+1},\o_{lcs},J_{can})$, $n\ge 1$,
admits a two-parameter deformation
by Vaisman lcK manifolds $(S^1\times S^{2n+1},\o_{lcs},J_{a,b})$, $b>0$. The group $T^2\times \SU (n+1) = S^1\times \U (n+1)$ acts
transitively on  $S^1\times S^{2n+1}$ preserving all of these lcK
structures. It is the maximal connected Lie group preserving any of the
above lcK structures. For $b\neq 1$ this group coincides with the
full connected isometry group of the lcK metric $g_{a,b}$.
For $b=1$ the full connected isometry group is strictly
larger, that is $\mathrm{Isom}_0(S^1\times S^{2n+1},g_{can}) = S^1\times
\SO (2n+2)$
\ec

\section{Homogeneous locally conformally symplectic manifolds}
Here we  give a  description of  homogeneous locally conformally
symplectic manifolds.

Let $(M=G/H,\o )$ be a homogeneous lcs
manifold with Lee form $\l$. For all of this section we will assume that $G$ is connected and effective and that 
$d\omega \neq 0$. We will consider $\o$ and $\l$ as
$\gh$-invariant forms on the Lie algebra $\gg$ which vanish on
$\gh$.

\subsection{A bound on the dimension of the center} 
\bp \label{boundProp} If $\l$ does not vanish on the center $\gz$ of $\gg$ then
$\dim \gz \le 2$. \ep

\pf As $\l$ is closed $\gg^\l:=\ker \l \subset \gg$ is an ideal.
Since $M$ is lcs we have the equation
$d\o = \l \wedge \o$ on $\gg$. Let $Z_0, Z_1\in \gz$, $\l (Z_0) =1$,
$Z_1, X\in \ker \l$. Then the above equation yields
\[ 0=d\o (Z_0,Z_1,X) = \o (Z_1,X).\]
This shows that $\gz \cap \gg^\l \subset \ker \o|_{\gg^\l}$, which implies
$\dim \gz \cap \gg^\l \le 1$ and, hence, $\dim \gz \le 2$.
\qed
\bc If $\gg$ admits an ad-invariant (possibly indefinite) scalar product
$b$ such that the vector $Z_0:=b^{-1}\l$ is not isotropic then
$\dim \gz \le 2$.
\ec
\pf It suffices to prove that $Z_0\in \gz$. For all $X,Y\in \gg$ we have:
\[ b([Z_0,X],Y)=b(Z_0,[X,Y])=\l ([X,Y]) = -d\l (X,Y) =0.\]
\qed

\bc If $G$ is reductive then $\dim Z(G) \le 2$. In particular,  a 
reductive  automorphism group of a homogeneous lcs manifold has at
most 2-dimensional center. \ec

\bp Let $(M=G/H,\o , g)$ be a homogeneous Vaisman manifold 
such that $G=\mathrm{Aut}(M,\o ,g)$. Then the center
$\gz$ of $\gg$ is $2$-dimensional.   
\ep 

\pf By Theorem \ref{VaismanThm}, the Reeb vector field is 
an infinitesimal automorphism of $(M,\o ,g)$, which generates
a one-parameter subgroup of $G$. Any vector $X\in \gg$ defines 
a Killing vector field $X^*$ on $M$. Let us denote by $Z\in \gg$ the 
{\cmssl Reeb vector}, that is the 
vector such that $Z^*$ is the Reeb vector field. Then the $G$-invariance
of $Z^*$ implies that $0=\mathcal{L}_{X^*}Z^*=[X^*,Z^*]=-[X,Z]$ for
all $X\in \gg$. Thus $Z\in \gz$, which implies $\dim \gz\ge 1$.   The same 
argument applies to the Lee field $\xi = -JZ$, showing that $\dim \gz \ge 2$. 
On the other hand,  Proposition \ref{boundProp} shows that $\dim \gz \le 2$. 
\qed  

 \subsection{A construction of  homogeneous lcs
manifolds }

Let $G$ be a Lie group  with  the Lie algebra $\gg$ and
 $Q = \mathrm{Ad}^*_G \phi = G/K$ the coadjoint orbit  of an element $\phi \in
 \gg^*$. We denote by $\omega_Q$ the (invariant) Kirillov-Kostant
 symplectic  form in $Q$  given by
 $$   (\omega_Q)_{\phi'}(X\cdot \phi',Y\cdot \phi') := \phi'([X,Y]), \,\, \phi' \in Q,\,\, X,Y \in \gg , 
 $$
where $X\cdot \phi' = -\phi' \circ ad_X\in T_{\phi'}Q$. Identifying $\o_Q$ with an $\mathrm{Ad}_K$-invariant $2$-form 
on $\gg$ vanishing on $\gk =\mathrm{Lie}\, K$ we can simply write 
\[ \o_Q(X,Y) = \phi ([X,Y]),\quad X,Y \in \gg .\]
 We will assume  that the orbit $Q$ is not conical, that is it is not
 invariant with respect to multiplication by positive numbers. Then
 the restriction $\phi|{\gk}$  of the  form $\phi$ to the
 stability subalgebra $\gk$ is not zero and $\gh := \gk \cap\ker \phi$
 is  an ideal of $\gk$ (see \cite{A}). We will assume that  the subalgebra $\gh$
 generates a closed  subgroup $H$ of $G$. Then we have:
 \bp (\cite{A}) The  $1$-form $\phi$  defines an invariant  contact  structure
 $\phi$ in $P = G/H$  and  the  contact manifold $(P = G/H,
 \phi)$ is a quantization of the homogeneous symplectic manifold
 $(Q = G/K, \omega_Q )$, that is $\phi$ is a connection on the 
$A$-principal bundle  $P = G/H \to G/K$ with the  curvature 
form $\omega_Q$, where $A = K/H \cong 
 \mathbb{R}$ or $\cong S^1$. 
\ep
 Let $D$ be a derivation of  the Lie algebra $\gg$  and $\gg(D) := \mathbb{R}D + \gg$
 the associated Lie  algebra  with  the ideal $\gg$. We denote by
 $\lambda$ the  closed $1$-form dual to $D$ (such that  $\lambda(D) =1, \lambda(\gg)
 =0)$ and define a $2$-form $\o$ on $\gg(D)$ by
 \be \label{fundEqu} \omega = -\lambda \wedge \phi + d
 \phi .\ee 
It is an $\ad^*_{\mathfrak{h}}$-invariant  $2$-form 
with kernel $\gh$ and satisfies
 $$    d \omega = \lambda \wedge d \phi = \lambda \wedge \omega .$$
We denote  by
 $G(D)$ a Lie group with the Lie algebra $\gg(D)$  and by $H$ its
 closed  (connected) subgroup generated by $\gh$. Obviously, we have:
 \bp  The $\mathrm{Ad}^*_H$-invariant  $2$-form $\omega$ defines an invariant
 lcs structure $\omega$ on the
 homogeneous manifold $M = G(D)/H$, that is an invariant non-degenerate $2$-form
 $\omega$  such that $d \omega = \lambda \wedge \omega$.
 \ep

We say that $(M = G(D)/H, \omega)$ is a homogeneous locally
conformally symplectic manifold  associated with the non-conical
orbit $Q = \mathrm{Ad}^*_G \phi$  and a derivation $D $ of the Lie algebra
$\gg$.

\noindent
{\bf Remark:} Let $(M,\o,J)$ be an lcK manifold of Vaisman type with 
Lee form $\l$ and Reeb form $\theta$. Then the equation \re{fundEqu} holds with 
$\phi=\frac{1}{2}\theta$. 

\subsection{The main result for homogeneous lcs manifolds}

In this subsection we show as a main result (Theorem \ref{mainlcsThm})
that the above construction gives all
homogeneous lcs manifolds satisfying a certain cohomological 
assumption, which we will explain now. 

Let $(M=G/H,\o )$ be a homogeneous lcs manifold with Lee form $\l$.
We consider $\omega$ and $\lambda$ as  $\mathrm{Ad}^*_H$-invariant forms on the Lie
algebra $\gg$, which vanish on $\mathfrak{h}$. 
Then $\o$ defines a cohomology class 
\[ {[}\o {]}\in H^2_\l (\gg ,\gh) :=\frac{\ker \left( 
d_\l : C^2(\gg ,\gh) \ra C^3(\gg ,\gh) \right) }{\mathrm{im} \left(  
d_\l:  C^1(\gg ,\gh)\ra  C^2(\gg ,\gh) \right)},\]  
where 
\[ C^k(\gg ,\gh) := \{ \a \in (\wedge^k\gg^*)^H| 
\iota_X\a =0\quad\mbox{for all}\quad X\in \gh\}\]
is the vector space of $\mathrm{Ad}^*_H$-invariant 
alternating $k$-forms vanishing on $\gh$ 
and \[ d_\l \a := d\a -\l\wedge\a ,\quad 
\mbox{for all}\quad \a \in \wedge^k\gg^*.\]  
We will assume that ${[}\o {]}=0$, which means that there exist
$\phi\in C^1(\gg ,\gh )$ satisfying the equation \re{fundEqu}. 
Recall that  $\gg':=
\mathfrak{g}^\l = \ker \lambda$ is an ideal of $\mathfrak{g}$ which
contains $\mathfrak{h}$. We can write
$$  \gg = \mathbb{R}D + \mathfrak{g}'$$
where $D \in \gg $  such that  $\lambda(D) =1$.
The assumption $d\omega \neq 0$ implies that $\lambda$ and $\phi$ are
linearly independent. Therefore, adding an element of $\gg'$ to
$D$, we can assume that $\phi (D)=0$.   
 The
restriction $\omega' = \omega|_{\mathfrak{g}'}$ is a closed $2$-form
on $\mathfrak{g}'$ and its kernel $\mathfrak{k}$ is a subalgebra
which contains  the codimension one   subalgebra $\mathfrak{h}$.
 \bl \label{fundEquLemma} Let $(M=G/H,\o )$ be a 
homogeneous lcs manifold with Lee form $\l$ and 
$d\o \neq 0$. Assume that $G$ contains the one-parameter subgroup
generated by the Reeb vector field $Z$ (see Proposition \ref{autoProp} and note that 
$Z$ is automatically complete since it is $G$-invariant). 
If ${[}\o {]}=0$ in $H^2_\l (\gg ,\gh)$ then 
the form $\omega$ can be written as
 \[  \omega = -\lambda \wedge \phi + d \phi , 
 \] 
where $\phi $ is an $\mathrm{Ad}^*_H$-invariant $1$-form on $\gg$ with $\ker
\phi \supset \mathbb{R}D + \mathfrak{h}$ which is not zero on
$\mathfrak{k}$. Moreover, 
\[ \o (Z,\cdot ) = \phi(Z) \l .\]
\el 

\pf 
Since ${[}\o {]}=0$, the equation 
\re{fundEqu} holds  for  some $\mathrm{Ad}^*_H$-invariant 
$1$-form $\phi$ which vanishes on $\gh$. The inclusion $\ker
\phi \supset \mathbb{R}D + \mathfrak{h}$ holds by our choice of $D$, 
as explained above.  We prove that $\phi|_{\gk}\neq 0$. 
Let  $Z\in \gg$ be the central element which corresponds to the Reeb vector field. 
Then $ad_Z^*\psi=0$ for every  $k$-form $\psi$ on $\gg$ and, in particular,  
\be   \iota_Zd\phi =  - ad_Z^*\phi=0 \label{adZthetaEqu}.\ee
Next we observe that the definition of the 
Reeb vector field (see Definition \ref{lcsDef}) implies that 
\be \label{lambda(Z)Equ}\l (Z)=0,\ee  
since $\o$ is skew-symmetric. Therefore  
the equations \re{fundEqu} and \re{adZthetaEqu} show that 
\be \label{omegaZEqu}\o(Z,\cdot ) = \phi(Z)\l .\ee
Since $\o$ is non-degenerate on $\gg/\gh$ this implies that 
\be \label{thetaZnot0}\phi (Z)\neq 0\ee 
and, hence, $\o (D,Z) = -\phi (Z)\neq 0$. 
So the plane $E$ spanned by $D$ and $Z$ is $\o$-non-degenerate.
Let $\gm'\subset \gg'$ be a subspace such that $\gm'\cap \gh =0$ 
and which projects to the $\o$-orthogonal complement of $\bar{E}=
(E+\gh)/\gh \subset \gg/\gh$ 
in $\gg/\gh$. In particular $\gm'\perp_\o Z$ implies 
\be \gg' =  \ker \lambda = \gh + \bR Z + \gm',\ee
in view of \re{lambda(Z)Equ} and \re{omegaZEqu}. 
Now we see that 
\be \gk = \ker \o' = \gh +\bR Z ,\ee
which, by \re{thetaZnot0}, proves that $\phi$ does not vanish on $\gk$.  
 \qed

We claim that the kernel $\gk$ of  the exact $2$-form
$\omega' = \o|_{\gg'}= d (\phi|_{\gg'})$ on $\gg'$  coincides  with the   
stabilizer of $\phi':=\phi|_{\gg'}$   in the  
coadjoint   representation of $\gg'$. 
In fact, this is a consequence of the equation
\[ \o' (X,\cdot ) = -\phi \circ ad_X|_{\gg'},\]
which holds for all $X\in \gg'$, in view of \re{fundEqu}. 
Hence,
the  corresponding  subgroup $K$ of the  group $G' \subset G$ is
closed.  By Lemma \ref{fundEquLemma}, the  
coadjoint orbit $Q := \mathrm{Ad}_{G'}^* \phi' = G'/K$ is
not conical  and $\gh = \gk \cap \ker \phi$ generates a closed
subgroup $H \subset G' \subset G$.  The $\mathrm{Ad}_H^*$-invariant $1$-form
$\phi'$ on $\gg'$ defines a  contact form  on $ P = G'/H$ and the
contact manifold $P = G'/H$ is a quantization of the  symplectic
manifold $Q = G'/K$. The contact property follows from the 
fact that  $d\phi'=\o'$ induces a non-degenerate $2$-form on 
$\gg'/\gk$ (see Lemma \ref{fundEquLemma}, and the next lemma).

\bl Under the assumptions of Lemma \ref{fundEquLemma}, we have 
\be \label{kerthetaEqu} \ker \phi' + \gk = \gg' .\ee
\el 

\pf  Since $\phi$ and $\l$ are linearly independent, 
$\phi' =\phi|_{\gg'} \neq 0$ and $\ker \phi' \subset \gg'$ is a
hyperplane. By \re{thetaZnot0}, 
$Z\not\in \ker \phi'$. Therefore, $\ker \phi' + \bR Z=\gg'$, 
which implies \re{kerthetaEqu}.  
\qed

Since $\ad_D|{\mathfrak{g}'}$ is a derivation
of the Lie algebra $\gg'$, we  can write $\gg = \gg'(\ad_D)$ and the
$2$-form $\omega$ on $\gg$ has the form
$$  \omega = - \lambda \wedge \phi + d \phi,$$
where $\phi$ is the canonical extension of $\phi'$ to 
a $1$-form on $\gg$. 
This shows:   
 \bt  \label{mainlcsThm} Any  homogeneous lcs
manifold satisfying the assumptions of Lemma \ref{fundEquLemma} can be  obtained  by the above construction,  that is  it
is associated  with   a non-conical  coadjoint orbit  $Q =
\mathrm{Ad}^*_{{G'}} \phi = G'/K$ of a Lie group $G'$ with the standard
symplectic form $\omega_Q = d \phi$   and a  derivation $D$ of the
Lie algebra $\gg'$. More precisely, it has the form $ ( M = G'(D)/H,
\omega )$ where  the Lie  algebra of $G'(D)$ is  the $D$-extension
$\gg'(D) = \mathbb{R}D + \gg'$ of  $\gg '$, $\gh := \ker \phi \cap
\gk $    and $\omega = -\lambda \wedge \phi + d \phi$.
 \et

Now we give some sufficient conditions which ensure the cohomological 
assumption used in this section.
\bd A homogeneous lcs manifold with Lee form $\l$
is called {\cmssl locally splittable} if the ideal $\gg'=\gg^\l\subset \gg$ has a complementary ideal, that is
$\gg = \bR D \oplus \gg' \;(D \in \gg)$.
It is called {\cmssl splittable} 
if $G=A\times G^\l$, where $A=\bR$ or $A=S^1$. 
\ed 
\bp \label{cohProp}Let $(M=G/H,\o )$ be a locally splittable 
homogeneous lcs manifold with Lee form $\l$ and 
$d\o \neq 0$. Then ${[}\o {]}=0$ in $H^2_\l (\gg ,\gh)$, 
$H^1_\l (\gg ,\gh)=0$ and $\dim Z(\gg' ) \le 1$. 
In particular, 
this is the case if
$\gg$ is reductive. 
\ep   

\pf We may assume that $\l (D) =1$. Then we decompose
$\o$ as
\be \label{decompEqu}\o = -\l\wedge \phi + \o',\ee
where $\phi$ and $\o'$ are $\mathrm{Ad}_H^*$-invariant forms on $\gg'$, 
which vanish on $\gh$. 
Differentiating this equation and comparing with the lcs equation, 
we obtain
\[ d\o = \l \wedge d\phi + d\o' = \l \wedge \o = \l\wedge \o'.\]
This shows that 
\[ \o'= d\phi .\]
Substituting this into \re{decompEqu} we get $d_\l \phi = \o$. 
To prove $H^1_\l (\gg ,\gh)=0$, let $\a \in C^1(\gg ,\gh )$ be a 
$d_\l$-closed form. We decompose it as
\[ \a = c\l + \a',\]
where $c$ is a constant and $\a'\in C^1(\gg' ,\gh ) \subset C^1(\gg , \gh )$. 
Differentiation yields
\[ 0=d_\l \a = -\l \wedge \a' + d\a',\]
which implies $\a'=0$ and $\a = c\l = - cd_\l 1$, where 
$1\in C^0(\gg ,\gh )=\bR$.   The bound on the dimension of the 
center of $\gg'$ follows from Proposition \ref{boundProp}. 
\qed

\bc Let $Q=G/K = \mathrm{Ad}_G^*\phi$ be a non-conical coadjoint orbit
such that the normal subgroup $H\subset K$ generated by 
$\gh = \ker \phi|_{\gk}$ is closed. Then $(P=G/H,\phi )$ 
is a homogeneous contact manifold and $(M=A \times P, 
\o = -dt\wedge \phi + d\phi )$ is a homogeneous 
lcs manifold, where $A=\bR$ or $A=S^1$. 
Conversely, any splittable homogeneous proper lcs manifold 
$(M=G/H,\o )$ with Lee form $\l$ 
can be obtained from this construction.  
\ec
We remark that the covering $\bR \times P$ of the lcs manifold 
$A\times P$  in the previous corollary, where $\bR\ra A$  is the 
universal covering group, is   
globally conformal to the symplectic cone over the contact manifold
$(P,\phi )$ after a redefinition $t=-2\tilde{t}$:\\ 
$\o = 2(d\tilde{t}\wedge \phi + \frac{1}{2}d\phi )
=\frac{2}{r^2}(rdr\wedge \phi +
\frac{r^2}{2}d\phi)$, where $\tilde{t}=\ln r$.  

\section{Homogeneous  locally conformally K\"ahler
manifolds of reductive groups}

\subsection{Left-invariant lcK structures on 4-dimensional reductive groups} 
\label{ExamplesSec} 
In this section we prepare the classification of  homogeneous lcK manifolds of reductive groups, 
to be given in Theorem \ref{lcKThm},  
by classifying left-invariant lcK structures 
on 4-dimensional reductive groups. 
We first describe all  left-invariant 
complex structures $J$ on such groups, then all left-invariant lcs structures
$\o$ and finally all left-invariant locally conformally pseudo-K\"ahler structures $(\o , J)$.
In particular, we describe all  lcK and Vaisman examples. This extends the results of \cite[Sec.\ 4]{HK2}. 
The following lemma  is a well known basic fact. 
\bl For any Lie group $G$,  the map 
\[ J\mapsto \gl_J:=\mathrm{Eig}(J,i) = \mathrm{ker} (J-i\mathrm{Id})\] 
induces a one-to-one correspondence between left-invariant complex
structures $J$ on $G$ and (complex)  Lie subalgebras $\gl =  \gl_J \subset \gg^\bC$ 
such that 
\be \label{lEqu} \gg^\bC = \gl+\rho \gl,\quad  \gl\cap \rho \gl=0,\ee
where $\rho$ denotes the real structure (i.e.\ complex anti-linear involutive automorphism) on $\gg^\bC$ with the fixed point
set $\gg$.
\el 

Let $\gg$ be a $4$-dimensional non-commutative reductive Lie algebra, that
is $\gg = \mathfrak{u} (2)$ or $\gg =\mathfrak{gl}(2,\bR )$, 
and $G$ any connected Lie group such that $\gg =\mathrm{Lie}\, G$. 
We may take $G=\mathrm{U}(2)$ or $G=\GL (2,\bR )$.  Let us denote by $\gg = \gz \oplus \gs$ the decomposition of 
the reductive Lie algebra $\gg$ into its center $\gz =\bR e_0$ 
and its maximal semisimple  ideal $\gs = [\gg , \gg ]$, which is $\su (2)$ or $\mathfrak{sl} (2,\bR )$. 
 We denote by $e^0$ the $1$-form on $\gg$ which 
vanishes on $\gs$ and has the value $e^0(e_0)=1$. 
\bl \label{BorelLemma} Let $G$ be a (connected) $4$-dimensional non-commutative reductive Lie group. 
Up to conjugation by an element of $G$, every left-invariant complex structure $J$ on $G$ is 
defined by a subalgebra $\gl_J = \mathrm{span}\{ e_0+e',e''\}$ such that
$e', e''\in \gs^\bC$, $[e',e''] = \mu e''$, $\mu\in \bC^*$. In particular, 
$e''$ belongs to the cone $\mathcal{C}\subset \mathfrak{sl}(2,\bC)$ of nilpotent elements. 
This is precisely the null cone with respect to the Killing form of $\mathfrak{sl}(2,\bC)\cong\bC^3$.
\el 
\pf 
We have to describe all subalgebras $\gl\subset \gg^\bC=\bC \oplus \mathfrak{sl}(2,\bC )$ satisfying \re{lEqu}. 
From $\r \gs^\bC=\gs^\bC$ we see that $\gl \not\subset \gs^\bC= \mathfrak{sl}(2,\bC)$.
Therefore $\gl$ admits a basis of the form $(e_0+e',e'')$, where $e',e''\in  \gs^\bC$. 
Then
\[ [e_0+e',e''] = [e',e''] \in \gl \cap  \gs^\bC = \bC e''\]
shows that 
\be \label{muEqu} [e',e''] = \mu e'',\quad \mu \in \bC^*.\ee 
Therefore $\mathrm{span}\{ e',e''\} \subset \gs^\bC$ is a Borel subalgebra and $e''$ 
belongs to the cone $\mathcal{C}$. \qed  
\bl \label{isotrLemma} Given a complex structure $J$ on $\gg$ and a $1$-form $\phi\in \gs^*\subset \gg^*$ such that 
$\o = e^0\wedge \phi +d\phi$ is non-degenerate (and, hence, defines
a lcs structure), the structure  $(\o , J)$ is 
locally conformally pseudo-K\"ahler if and only if $\gl_J= \mathrm{span}\{ e_0+e',e''\}\subset \gg^\bC$ 
is isotropic with respect to $\o$. This is the case if and only if either 
$\mu =1$ or $\phi (e'')=0$. 
\el 
\pf
Notice first that the $2$-form $\o$ is $J$-invariant if and only if
it is of type $(1,1)$, which means that $\gl_J$ and $\r \gl_J$ are
isotropic. Next we evaluate $\o = e^0\wedge \phi +d\phi$ on the basis of $\gl_J$: 
\[ \o (e_0+e',e'')= \phi (e'') -\phi ([e',e'']) = (1-\mu )\phi (e'').\]
\qed 
\subsubsection*{The compact case}
Let us first consider the case $\gs=\su (2)$  and denote by $(e_1,e_2,e_3)$ a basis of $\su (2)$ such that
$[e_\a,e_\b]=-e_\g$ for every cyclic permutation of $(1,2,3)$. In the following 
$(\a ,\b,\g)$ will be always a cyclic permutation. 
Then the basis $(e^0,e^1,e^2,e^3)$ of $\gg^*=\mathfrak{u}(2)^*$ which is dual to
$(e_0,e_1,e_2,e_3)$ has the following differentials:
\[ de^0 =0,\quad de^\a = e^{\b\g}:= e^{\b}\wedge e^\g.\]

\bp Up to conjugation by an element of $\U (2)$, every left-invariant complex structure $J$ on $\U (2)$ is contained
in the following Calabi-Eckmann family
\be \label{CEEqu} Je_0 = ae_0 +be_1, Je_1 = ce_0-ae_1,Je_2=-e_3,Je_3=e_2,\ee
which depends on two-parameters $a\in \bR$ and $b\neq \bR^*$; $c=-\frac{1+a^2}{b}$. 
\ep 
\pf We specialize the description of complex structures in Lemma \ref{BorelLemma}. 
Since $\U(2)$ acts transitively on the quadric $Q=P(\mathcal{C})\cong \bC P^1$ we can assume that
$e''=e_2+ie_3$. Then the equation \re{muEqu} shows that $e'\equiv -i\mu e_1\pmod{\bC e''}$
and we can choose the above basis of $\gl$ such that $e'=-i\mu e_1$.  Then \re{lEqu}
is satisfied if and only if $\r e'\neq e'$, i.e.\ $\mu\not\in i\bR$. This shows that the complex
structure $J$ defined by $\gl_J=\gl$ is given by \re{CEEqu}, where $\mu = \mu_1 +i\mu_2$ is 
related to $a,b,c$ by 
\be  \label{amuEqu} a= \frac{\mu_2}{\mu_1},\quad b= \frac{|\mu|^2}{\mu_1},\quad c= -\frac{1}{\mu_1}.\ee
\qed 
\bp Up to scale, every left-invariant lcs form on $\U(2)$ is of the 
form
\be \label{e0EquU2}\o = e^0\wedge \phi +d\phi,\ee
where $\phi= \sum a_\a e^\a\in \gs^*$ is any nonzero form. All these structures
are equivalent up to conjugation in $\U (2)$. 
\ep 

\pf Let $\o$ be an lcs structure on $\gg=\mathfrak{u}(2)$. Since $e^0$ is the only closed $1$-form on $\gg$, up to scale, 
we can assume that the 
Lee form of $\o$ is given by $\l = -e^0$. The canonical $1$-form of $\o$ is 
given by a nonzero element $\phi\in \gs^*$ and any such element defines
an lcs structure $\o$ by the formula \re{e0EquU2}. \qed 

\bt \label{u2Thm} Let $J=J_{a,b}$ be any of the left-invariant complex structures on $G=\U(2)$, as defined in \re{CEEqu}. 
\begin{enumerate}
\item[(i)]
If $(a,b)\neq (0,1)$ then, up to scale, there is a unique left-invariant lcs structure $\o$ on  $\U(2)$ such that
$(\o ,J)$ is locally conformally pseudo-K\"ahler. It is given by $\o = e^{01} +e^{23}$. All these structures
are of Vaisman type. The locally conformally pseudo-K\"ahler metric $g=-\o \circ J$ is definite if  and only if $b<0$. 
\item[(ii)]If $(a,b)=(0,1)$ then $(\o ,J)$ is locally conformally pseudo-K\"ahler for 
every left-invariant lcs structure $\o$ on  $\U(2)$. The metric is always indefinite and the structure $(\o ,J)$ is of Vaisman
type if and only if $\o$ is proportional to  $e^{01}+e^{23}$. 
\end{enumerate}
\et 
\pf The pair $(\o , J)$ defines a locally conformally pseudo-K\"ahler structure on $G$ if and only if
$\gl_J= \mathrm{span}\{ e_0+e',e''\}\subset \gg^\bC$ is isotropic with respect to $\o$, where $e'=-i\mu e_1$,  
$e''= e_2+ie_3$. To check
this property we evaluate \re{e0EquU2},
\be \label{oabcEqu} \o = -\l \wedge \phi + d\phi = \sum a_\a e^{0\a} +\sum a_\a e^{\b\g}\ee
on the above basis of  $\gl_J$: 
\begin{eqnarray*} \o (e_0+e',e'') &=& a_2  + ia_3 + a_2e^{31}(-i\mu e_1,ie_3) +a_3e^{12}(-i\mu e_1,e_2)=
a_2  + ia_3  -\mu a_2-i\mu a_3\\
&=& (1-\mu)(a_2+ia_3 ).\end{eqnarray*}
So we see that $\gl_J$ is $\o$-isotropic if and only if
either 
\begin{enumerate}
\item[(i)] $a_2=a_3=0$, that is $\o= e^{01} +e^{23}$, up to scale, or
\item[(ii)]  $\mu=1$, that is   $(a,b)=(0,1)$.  
\end{enumerate}
In case (i) we compute
\[ 2\xi = \o^{-1}J^*\l = -\o^{-1} (ae^0+ce^1) = -(-ae_1+ce_0)= ae_1-ce_0\]
and 
\[ 2Z = 2J\xi = a(ce_0-ae_1) -c(ae_0+be_1)= (-a^2-cb)e_1=e_1.\]
This shows that $X=2(\xi -aZ)=-ce_0\in \gz$ and, hence, defines a (nonzero) Killing vector field. 
On the other hand, $\mathcal{L}_v\o=0$ for all $v\in \mathrm{span}\{ e_0,e_1\} = \mathrm{span}\{Z,\xi \}$,
since $e_0,e_1\in \ker d\phi = e^{23}$, where 
\[  \mathcal{L}_v := d \circ \iota_v + \iota_v\circ d : \wedge^k \gg^* \ra   \wedge^k \gg^*\]  
is the linear map induced by the Lie derivative in direction of the left-invariant vector
field associated with the vector $v\in \gg$. 
In particular, $\mathcal{L}_X\o=0$. 
These two properties of $X$ show that $X$ and, therefore, $JX$ define (real)  holomorphic vector fields.
Writing $\xi$ as a linear combination of $X$ and $JX$ we see that also 
$\xi$ defines a holomorphic vector field. On the other hand,  by the same argument as for $X$ we see that 
$\mathcal{L}_\xi \o=0$, since 
$\xi$ is a linear combination of $e_0$ and $e_1$. Therefore $\xi$ defines a Killing vector field. 
Now it suffices to remark that a locally conformally pseudo-K\"ahler manifold is Vaisman if and only if the Lee field 
is Killing. In fact, the Lee field is locally a gradient vector field (due to $d\l=0$) and  a gradient vector field
is Killing if and only if it is parallel. To finish the proof of (i) we have 
to check when the metric $g=-\o \circ J$ is definite. We compute
\begin{eqnarray*}  \o \circ J&=& J^*e^0\ot e^1 - J^*e^1\ot e_0 + J^*e^2\ot e^3-J^*e^3\ot e^2\\
&=&
(ae^0+ce^1)\ot e^1 -(be^0-ae^1) \ot e^0 + e^3\ot e^3 +e^2\ot e^2\\
&=& -b(e^0)^2+2ae^0e^1+c(e^1)^2 + (e^3)^2 +(e^2)^2,\end{eqnarray*}
which is definite if and only if $b<0$. 
To prove (ii) we compute $\o \circ J$ for $\o$ given in 
\re{oabcEqu} and $J=J_{0,1}$:
\begin{eqnarray*}  \o \circ J&=& \sum a_\a (J^*e^0\ot e^\a -J^*e^\a \ot e^0) + 
\sum a_\a (J^*e^\b \ot e^\g -J^*e^\g \ot e^\b )\\
&=& -\sum a_\a e^1\ot e^\a -a_1(e^0)^2-a_2e^3\ot e^0 +a_3e^2\ot e^0 +
a_1((e^2)^2 +(e^3)^2)\\
&& -a_2(e^2\ot e^1 + e^0\ot e^3)+a_3(e^0\ot e^2 - e^3\ot e^1)\\  
&=& -a_1(e^1)^2 -a_1 (e^0)^2+a_1(e^2)^2 +a_1(e^3)^2 -2a_2e^1e^2-2a_3e^1e^3-2a_2e^3 e^0 +2a_3e^2e^0.
\end{eqnarray*}  
This metric is always of signature $(2,2)$. 
Now suppose that $(\o , J)$ is of Vaisman type. Then the Lee vector $\xi$ satisfies
$\mathcal{L}_\xi \phi=\iota_\xi d\phi=0$. This implies that $\xi$ is a linear 
combination $c_0e_0+c_1\vec{a}$ of $e_0$ and 
$\vec{a}=\sum a_\a e_\a$. Since $g(\xi, \cdot ) = -\frac12\l$ applying 
$\o \circ J$ to $c_0e_0+c_1\vec{a}$ should be a multiple of $\l=-e^0$. 
We calculate
\begin{eqnarray*} \o J(c_0e_0+c_1\vec{a}) &=& c_0(-a_1e^0-a_2e^3+a_3e^2) +
c_1 a_1(-a_1e^1-a_2e^2-a_3e^3)\\ 
&&+c_1a_2(a_1e^2-a_2e^1+a_3e^0)+
c_1a_3(a_1e^3-a_3e^1-a_2e^0).\end{eqnarray*}
The coefficient of $e^1$ is 
\[ -c_1\sum a_\a^2\]
and has to vanish. Since $\vec{a}\neq 0$ this shows that $c_1=0$ and that 
$\xi$ is proportional to $e_0$.  Then 
\[  \o Je_0 = -a_1e^0-a_2e^3+a_3e^2,\]
which is proportional to $e^0$ only if $a_2=a_3=0$. 
This implies $\o = e^{01}+e^{23}$ up to a factor, as claimed. 
\qed 
\subsubsection*{The non-compact case}
Let us now consider the case $\gs=\mathfrak{sl} (2,\bR )$  and denote by $(h,e_+,e_-)$ a basis of $\mathfrak{sl}(2,\bR )$ such that
$[h,e_\pm]=\pm 2e_\pm$, $[e_+,e_-]=h$. Then the basis $(e^0,h^*,e^+,e^-)$ of $\gg^*=\mathfrak{gl}(2,\bR )^*$ 
which is dual to  $(e_0,h,e_+,e_-)$ has the following differentials:
\[ de^0=0,\quad dh^*=-e^+\wedge e^-,\quad de^\pm =\mp 2h^*\wedge e^\pm .\]
We denote by $\r$ the standard real structure on $\gg^\bC$ associated with 
the real form $\gg=\mathfrak{gl}(2,\bR )$. 
\bp \label{GL2cxProp} Up to conjugation by an element of $\GL (2,\bR )$, every left-invariant complex structure $J$ on $\GL (2,\bR )$ 
belongs to one of the following two families depending on $\mu=\mu_1+ i \mu_2\in \bC \setminus i\bR$. 
\begin{enumerate}
\item[(i)] 
\begin{eqnarray*} Je_0&=& \frac{\mu_2}{\mu_1}e_0 -\frac{|\mu|^2}{2\mu_1}(e_+-e_-)\\
Jh &=&e_++e_-\\
Je_\pm &=& \pm \frac{1}{\mu_1}e_0 \mp \frac{\mu_2}{2\mu_1}(e_+-e_-)-\frac12 h.
\end{eqnarray*} 
\item[(ii)] 
\begin{eqnarray*} Je_0&=& \frac{\mu_2}{\mu_1}e_0 +\frac{|\mu|^2}{2\mu_1}(e_+-e_-)\\
Jh &=&-(e_++e_-)\\
Je_\pm &=& \mp \frac{1}{\mu_1}e_0 \mp \frac{\mu_2}{2\mu_1}(e_+-e_-)+\frac12 h.
\end{eqnarray*} 
\end{enumerate}
These two families are related by the outer automorphism
of $\mathfrak{gl}(2,\bR )$ which maps $(e_0,h,e_\pm )$ to $(e_0,h,-e_\pm )$. (See 
remark below for a description of these complex structures in a  basis which is orthonormal with respect to a 
suitably normalized bi-invariant scalar product on $\mathfrak{gl}(2,\bR )$.) 
\ep 
\pf As before, any complex structure is defined by a subalgebra $\gl\subset \gg^\bC$ satisfying \re{lEqu}.
The latter admits a  basis $(e_0+e',e'')$, where $e',e''\in  \gs^\bC$.  Then $[e',e''] = \mu e''$,  $\mu \in \bC^*$,
and $e''\in \mathcal{C}$. The group $\SL (2,\bR )$ has three orbits on the 
quadric $Q=P(\mathcal{C})$.  As representatives $e''$ of these orbits we choose
\[ e_+,\quad ih + e_++e_-,\quad h +i(e_++e_-).\]
The first case is excluded, since $\r e_+ = e_+$. 
The elements $e'$ corresponding to $e''= ih + e_++e_-$ and 
$e''=h +i(e_++e_-)$
are given by
\[  \frac{i\mu}{2}(e_+-e_-),\quad - \frac{i\mu}{2}(e_+-e_-).\]
Again $\mu \not\in i\bR$ by  \re{lEqu}. This gives the two families (i) and (ii). 
\qed 
Using the Killing form we can identify $\gs^*$ with $\gs$. Since the Killing form of $\gs = \mathfrak{sl}(2,\bR )$  is 
Lorentzian we can further identify $\gs$ with a Lorentzian vector space $\bR^{2,1}$. 

\noindent
{\bf Remark:} Putting $e_1:= (e_+-e_-)/2$,  $e_2=h/2$, $e_3:= (e_++e_-)/2$ and using the abbreviations \re{amuEqu} we can rewrite 
the complex structures in Proposition \ref{GL2cxProp} in a form similar to \re{CEEqu}:
\begin{enumerate}
\item[(i)] 
\[ Je_0 = ae_0 -be_1,\quad
Je_1 = -ce_0-ae_1,\quad
Je_2 = e_3,\quad 
Je_3 = -e_2.  
\] 
\item[(ii)] 
\[ Je_0 = ae_0 +be_1,\quad
Je_1 = ce_0-ae_1,\quad
Je_2 = -e_3,\quad 
Je_3 = e_2.  
\] 
\end{enumerate}

\bp \label{lcsGL2Prop} Up to scale, every left-invariant lcs form on $\GL(2,\bR)$ is of the 
form
\be \label{e0Equ}\o = e^0\wedge \phi +d\phi,\ee
where $\phi= \sum a_\a e^\a\in \gs^*\cong \gs = \mathfrak{sl}(2,\bR )= \bR^{2,1}$ is any non-isotropic $1$-form.  \ep 
\pf
It suffices to check that $\o$ is non-degenerate if  and only if $\phi$ is space-like or time-like.
\qed 

Next we describe all left-invariant lcs structures which are
compatible with any of the complex structures $J_\mu$ on $G=\GL(2,\bR)$, as described in 
Proposition \ref{GL2cxProp}. It is sufficient to consider the family (i), since  it is equivalent to 
(ii) by an automorphism of $G$. 

\bt \label{GL2Thm} Let $J=J_\mu$ be any of the left-invariant complex structures on $G=\GL(2,\bR)$, as defined in 
Proposition \ref{GL2cxProp} (i). 
\begin{enumerate}
\item[(i)]
If $\mu\neq 1$ then, up to scale,  there is a unique left-invariant lcs structure $\o$ on 
$\GL (2,\bR)$ such that
$(\o ,J)$ is locally conformally pseudo-K\"ahler. It is given by\linebreak
 \[ \o =e^0\wedge(e^+-e^-) -2h^*\wedge (e^++e^-)= e^0\wedge e^1-e^2\wedge e^3, \]  
 where $(e^0,e^1,e^2,e^3)$ denotes the basis dual to $(e_0,e_1,e_2,e_3)$. 
 All these structures
are of Vaisman type with (positive or negative) definite metric. 
\item[(ii)] If $\mu=1$ then $(\o ,J)$ is locally conformally pseudo-K\"ahler for 
every left-invariant lcs structure $\o=e^0\wedge\phi +d\phi$ on  $\GL(2,\bR )$. 
The locally conformally pseudo-K\"ahler metric $g= -\o \circ J$ associated with a non-isotropic $1$-form 
$\phi = a_hh^*+a_+e^++a_-e^-\in \gs^*$
is given by
\begin{eqnarray} g &=&  -\frac12(a_+-a_-)(e^0)^2-2(a_+-a_-)(h^*)^2+2(a_++a_-)e^0h^* 
-2a_+(e^+)^2 +2a_-(e^-)^2\nonumber \\
&&- a_he^0(e^++e^-) -2a_hh^*(e^+-e^-).\label{metricEqu}\end{eqnarray}
It is of Vaisman type if and only if $a_h=0$ and $a_+=-a_-\neq0$, in which case the metric is definite.
In particular, the  locally conformally pseudo-K\"ahler metric $g$ is  non-Vaisman and positive definite  
if and only if, first, $a_h\neq 0$ or 
$a_+\neq -a_-$ and, second, $- a_h^2 > 4\, a_+ a_-$ and $a_- >0> a_+$.
\end{enumerate}
\et 
\pf 
According to Proposition  \ref{lcsGL2Prop} any lcs structure on $\gg$ is of the form
$\o = e^0\wedge \phi +d\phi$, where $\phi= a_hh^*+a_+e^++a_-e^-\in \gs^*$ is any non-isotropic $1$-form.  
It is of type $(1,1)$ with respect  to $J$ if and only if either (i) $\phi (e'')=a_h +i(a_++a_-)=0$ or 
(ii) $\mu=1$ (see Lemma \ref{isotrLemma}).  In the first case, we have, up to scale, 
$\phi = e^+-e^-$, which implies $\o=e^0\wedge(e^+-e^-) -2h^*\wedge (e^++e^-)$.  
The corresponding locally conformally pseudo-K\"ahler metric 
$g$ is definite and Vaisman (the above basis
of $\gg$ is $g$-orthogonal). In the second case, a straightforward calculation of the metric 
yields  the above formula \re{metricEqu}, depending on the parameters 
$a_h,a_\pm$. Assuming that this metric is Vaisman, we see that 
\[\xi \in \ker d\phi= \mathrm{span}\{ e_0,\vec{a}=\frac{a_h}{2}h+a_+e_-+a_-e_+\}.\] 
So $\xi = \a e_0 +\b \vec{a}$ for some $(\a,\b)\in \bR^2\setminus \{ 0\}$.  Then using 
\re{metricEqu} we see that $g(\x ,\cdot )$ is proportional to $\l=-e^0$ if and only if
the following equations hold
\begin{eqnarray*}
\a (a_++a_-)=0\\
\a a_h=0\\
\b (\frac{a_h^2}{2}+2a_+a_-)=0.
\end{eqnarray*} 
Since $\phi$ is not light-like, we see that $\frac{a_h^2}{2}+2a_+a_-\neq 0$. Therefore
$\b=0$ and $\a\neq0$, which shows that $a_h=a_++a_-=0$. In that case, 
$g= -a_+(e^0)^2-4a_+(h^*)^2-2a_+(e^+)^2-2a_+(e^-)^2$, which is definite.  
Now it suffices to check that the metric \re{metricEqu} is always definite if 
$a_h=0$ and $a_+a_-<0$. (In the case $a_+<0$ it is positive definite.) 
Now that we have characterized the Vaisman case in (ii), it follows that 
the metric is non-Vaisman if and only if $a_h\neq 0$ or 
$a_+\neq -a_-$. So it only remains to check that the metric is positive definite
if and only if $- a_h^2 > 4\, a_+ a_-$ and $a_- > 0> a_+$.  
This is obtained from
a calculation of principal minors.\qed 

\subsection{Classification of  
homogeneous  lcK manifolds of reductive groups}

In this subsection we prove the following main theorem.
\bt  \label{lcKThm} Every homogeneous proper lcK manifold  $(M=G/H,\omega ,J)$ of 
a connected reductive Lie group $G$ such that $H$ is connected and $N_G(H)$ is compact is of Vaisman type.
\et

\pf We assume without restriction of generality that $G$ is effective. 
As before we consider the fundamental form $\o$, the Lee form $\l$ and the Reeb form $\theta = \frac{1}{2}J^*\l$ as $H$-invariant 
forms on $\gg$ which vanish on $\gh$. By Proposition \ref{cohProp} we know that there exist $\phi \in C^1(\gg,\gh)$ such that 
\re{fundEqu} is satisfied and that the $1$-form $\phi$ is unique up to addition of a multiple of $\l$. 
Let $\gm\subset \gg$ be an $H$-invariant complement of $\gh$ containing
the center $\gz$ of $\gg$. Let us denote by $Z, \xi\in \gm$ the linearly independent $H$-invariant vectors which correspond 
to the Reeb and Lee vector fields on $M$. We choose $\phi$ such that $\phi (\xi ) =0$. Together with 
the equation \re{fundEqu} this makes $\phi$ unique. We will call $\phi$ the {\cmssl canonical $1$-form}.

\bp \label{phiProp}Under the assumptions of  Theorem \ref{lcKThm}, the canonical $1$-form 
coincides with the Reeb form $\theta$ up to a factor $1/2$:
\[ \phi = \frac{1}{2}\theta.\] 
\ep 
\pf  
The proof of Proposition \ref{phiProp} is based on the following key lemma, the proof of which is given below. 
\bl \label{keyLemma} Under the assumptions of  Theorem \ref{lcKThm}, we have $Z, \xi \in \ker d\phi$.
\el 
Using Lemma \ref{keyLemma}, we compute
\[ \mathcal{L}_\xi \phi = \iota_\xi d\phi=0,\]
where, for any $\mathrm{Ad}_H$-invariant $v\in \gm$, 
\[  \mathcal{L}_v := d \circ \iota_v + \iota_v \circ d : C^k (\gg,\gh)  \ra   C^k (\gg ,\gh) .\] 
$\mathcal{L}_v$ is the 
linear map induced by the Lie derivative in direction of the $G$-invariant vector
field $X_v$ which extends $v$. 
Since also  $\mathcal{L}_\xi \l = \iota_\xi d\l =0$, the equation \re{fundEqu} implies 
\be\label{LxiEqu}  \mathcal{L}_\xi \o = - \l \wedge \mathcal{L}_\xi \phi + d\mathcal{L}_\xi \phi=0.\ee
Now Lemma \ref{LxioLemma} shows that 
\[ \o =  -\frac{1}{\l (\xi )}d_\l \theta = \frac{1}{2}d_\l  \theta .\]
Since $\o = d_\l \phi$ and $H^1_\l (\gg ,\gh )=0$, this proves that $\phi = \frac{1}{2}\theta \pmod{\bR \l}$.
Finally, for the canonical $1$-form we have $\phi (\xi )=0$, 
such that $\phi = \frac{1}{2}\theta$.
This finishes the proof of Proposition \ref{phiProp}. 
\qed 

\pf (of Lemma \ref{keyLemma}) 
Let us denote by $G_0$ the maximal connected subgroup
of the normalizer of $H$ in $G$. Since $H$ is compact,
$G_0$ is reductive. The Lie algebra $\gg_0$ of $G_0$ is decomposed
as 
\[ \gg_0 = \gh + \gm_0,\]
where $\gm_0=Z_\gm (\gh )$ contains $\gz$, $Z$ and $\xi$. 
Since $J$ is $H$-invariant, the maximal trivial $H$-submodule $\gm_0 \subset \gm$
is $J$-invariant. This implies that $\o$ is non-degenerate on $\gm_0$,  because
$g=-\o \circ J$ is positive definite.  Therefore the restriction of 
$(\o , J)$ to $\gm_0$ defines an invariant lcK structure on $M_0=G_0/H$
with the Lee form $\l_0 = {\l}|_{\gm_0}$.  Notice that $\l_0\neq 0$, since $\xi\in \gm_0$. 
Therefore, the lcK structure on $M_0$ is not K\"ahler, unless $\dim M_0=2$. 
From the fact that $H$ is normal in $G_0$, we see that $M_0$ is a Lie group. 
In the K\"ahler case, the Lie group $M_0$ is 2-dimensional and thus Abelian. 
So, in that case, $d\phi =0$ and the assertion of Lemma \ref{keyLemma} follows. 
Otherwise $M_0$ is at least 4-dimensional and the lcK structure is non-K\"ahler. Therefore, we can assume
from the beginning that $H$ is trivial.  This reduces the proof of 
Lemma \ref{keyLemma} to the following special case.

\bl  \label{keyLemmaspecial} Under the assumptions of  Theorem \ref{lcKThm} and
the additional assumption that $H$ is trivial, we have $Z, \xi \in \ker d\phi$.
\el 
\pf Let $B$ be a non-degenerate $\mathrm{Ad}_G$-invariant symmetric bilinear form on $\gg$. Then 
there exists endomorphisms $A_\o, A_g, A_{d\phi}, A_{\l \wedge \phi}\in \End\gg$ and a vector $v=v_\phi \in \gg$ such that  
\[ \o = B\circ A_\o,\quad g= B\circ A_g,\quad d\phi = B\circ A_{d\phi},\quad \l \wedge \phi=B\circ A_{\l \wedge \phi},\quad \phi = Bv. \]
We claim that  
\[ A_{d\phi}=-ad_v,\quad A_{\l \wedge \phi}= \l \ot v +2\phi \ot A_g\xi.\]
In fact,
\[ d\phi = -\phi \circ [\cdot ,\cdot ] = -B(v, [\cdot ,\cdot ] ) = B([\cdot ,v],\cdot ) = -B\circ ad_v,\]
\[ \l \wedge \phi = \l \ot \phi -\phi \ot \l = \l \ot Bv -
\phi \ot (-2g\xi)= B\circ (\l \ot v + 2\phi \ot A_g\xi).\]
 The equation $\o = -\l \wedge \phi + d\phi$ can now be rewritten as
 \[ A_\o = -A_{\l \wedge \phi} -ad_v =-\l \ot v -2\phi \ot A_g\xi -ad_v.\]
Since $\l$ and $\phi$ are linearly independent ($d\o\neq 0$), the skew-symmetric
endomorphism $A_{\l \wedge \phi}$ has rank two. More precisely, 
\[\mathrm{im}\, A_{\l \wedge \phi}=
\mathrm{span}\{ v,A_g\xi\}.\] 
Notice that $-2 (B \circ A_g) \xi = -2 g\xi = \l$. Therefore, the equation $d\l =0$ shows
that $A_g\xi \in \gz=[\gg ,\gg ]^{\perp_B}$.  In particular, $\gz \neq 0$. Since $A_\o$ has maximal rank, we see
that the image of $ad_v$ is complementary to $\mathrm{span}\{ v,A_g\xi\}$ in $\gg$ and
of codimension one in the semisimple Lie algebra $\gs = [ \gg ,\gg ]\supset \mathrm{im}\, ad_v$. 
This implies that the centralizer $Z_\gs (v)$ of $v$ in $\gs$  is one-dimensional. 

This shows that the rank of $\gs$ is one and $\dim \gs =3$. 
Since the dimension of $\gg$ is even,  the inequality $1\le \dim \gz\le 2$ implies that $\dim \gz=1$. 
Therefore,  $\gg = \mathfrak{u}(2)$, because $\gg$ is compact.  

We have proven in Section \ref{ExamplesSec} 
that all lcK structures on  $\gg = \mathfrak{u}(2)$ are of Vaisman type
and, hence,  satisfy $Z, \xi \in \ker d\phi$. 
This finishes the proof of Lemma \ref{keyLemmaspecial} and Lemma \ref{keyLemma},
and thus completes the proof of Proposition \ref{phiProp}. 
\qed 

The following Proposition finishes the proof of Theorem \ref{lcKThm}. 
\qed 

\bp Let $(M=G/H,\o , J)$ be a homogeneous proper lcK manifold of 
a reductive  Lie group $G$ such that $N_G(H)$ is compact and such that the canonical $1$-form is given by $\phi = \frac{1}{2}\theta$. 
Then $(M=G/H,\o , J)$ is of Vaisman type.
\ep  

\pf  Using  the assertion $\xi \in \ker d\phi$ in Lemma \ref{keyLemma}, we have shown in \re{LxiEqu} that $  \mathcal{L}_\xi \o =0$. Similarly, 
$Z\in \ker d\phi$ implies 
\[ \mathcal{L}_Z \phi = \iota_Z d\phi=0\]
and, hence, 
\[  \mathcal{L}_Z \o = -\l \wedge \mathcal{L}_Z \phi + d\mathcal{L}_Z \phi=0.\] 
We claim that 
\be \mathrm{span}\{ Z,\xi\}\cap \gz \neq 0.\ee
Since $Z$, $\xi$ and $\gz$ are contained in 
the normalizer $\gg_0=N_{\gg}(\gh)$ of $\gh$ in $\gg$, it is sufficient to prove this in the case  $\gg=\mathfrak{u}(2)$, $\gh=0$. 
Recall that any element $X\in \gg$ defines a Killing vector field $X^*$ on $M=G/H$ and that any 
 $\mathrm{Ad}_H$-invariant element $X\in \gm$ extends
as a $G$-invariant vector field $\tilde{X}$ on $M=G/H$.  
If $X\in \gz\subset \gm$ then $\tilde{X}=X^*$, that is 
$\mathcal{L}_{\tilde{X}}g=0$.  
If $0\neq X\in  \mathrm{span}\{ Z,\xi\}\cap \gz$, then $\mathcal{L}_{Z}\o=\mathcal{L}_{\xi}\o=0$ imply 
$\mathcal{L}_{X}\o =0$ and, hence, 
$\mathcal{L}_{\tilde{X}}\o =0$.   Combining these equations, we see that $\mathcal{L}_{\tilde{X}}J =0$, which implies
that the Reeb and the Lee vector fields are both holomorphic. Since the Lee field is a gradient vector 
field ($d\l =0$) this shows that the Lee field is parallel. This proves the proposition. \qed 

\noindent 
{\bf Example:} Note that the normalizer
 $N_G(H)=T^2=S^1\times S^1$ of $H=\SO(2)\subset \SL(2,\bR)$ in $T^2 \times \SL (2,\bR)$ is compact.
Therefore, Theorem \ref{lcKThm}  shows that every 
$G$-invariant lcK structure on $M=G/H=T^2\times
\SL(2,\bR)/\SO(2)$ is of Vaisman type. This should be contrasted with the fact that 
$S^1\times \SL (2,\bR )$ admits left-invariant non-Vaisman lcK structures by 
Theorem \ref{GL2Thm}.

 \subsection{Left-invariant lcK structures on reductive Lie groups}
 In this section we specialize to the case of  left-invariant lcK structures
 on Lie groups $G$. We will not assume that $G$ is compact and will allow the 
 pseudo-K\"ahler metric to be indefinite. 

 \bt \label{linvThm} Let $(G,\o ,J)$ be a Lie group endowed with a left-invariant (proper) 
 locally conformally pseudo-K\"ahler structure. 
 \begin{enumerate}
 \item[(i)] If $\gg = \mathrm{Lie}\, G$ admits a bi-invariant (possibly indefinite) scalar product $B$ with
non-isotropic $B^{-1} \lambda$, then 
 the dimension of the centralizer of $v$ (as defined in Lemma \ref{keyLemmaspecial}) 
 in $\gg$ is at most $2$. 
 \item[(ii)] If $\gg$ is reductive, then we have either $\gg = \mathfrak{u}(2)$ or $\gg = \mathfrak{gl}(2,\bR )$,
 and $(\o ,J)$ is one of the 
 locally conformally pseudo-K\"ahler structures classified in Theorems \ref{u2Thm}  and \ref{GL2Thm}. 
 In both cases  there exist  locally conformally pseudo-K\"ahler structures that are not
 of Vaisman type and in the case $\gg = \mathfrak{gl}(2,\bR )$ there even exist such structures
 that are not of Vaisman type with positive definite metric. 
 \end{enumerate}
 \et
 
 \pf We keep the same notation as in the proof of Lemma \ref{keyLemmaspecial}. 
 We first note that since $B^{-1} \lambda$ is non-isotropic, $\gg$ is splittable; and thus 
  ${[}\o {]}=0$ in $H^2_\l (\gg)$. The equation \[ ad_v=-A_\o -\l \ot v -2\phi \ot A_g\xi \]
 proven there (without using the compactness assumption of Lemma \ref{keyLemmaspecial}) shows that 
 the rank of $ad_v$ is at least $\mathrm{rk}\, \o -2 = \dim \gg -2$.  This implies that $Z_\gg (v)$ is at most two-dimensional.
 This proves (i).  Now we prove (ii). If $\gg$ is reductive the image of 
 $ad_v$ is necessarily a proper subspace of $\gs$. To see this it is sufficient
 to decompose $v$ according to the decomposition $\gg = \gs \oplus \gz$. 
 This proves that the image of $ad_v$ in $\gs$ is a hyperplane and that $Z_\gs (v)$ is one-dimensional, since 
 $0\neq A_g\xi \in \gz$.  Since the nilpotent part as well as the semisimple part of $ad_v|_\gs$ belongs
to $Z_\gs (v)\subset \gs \cong \ad (\gs)$, it follows that $ad_v|_\gs$ is either semisimple or nilpotent. 
It is clear that the dimension of the centralizer of a semisimple element in a semisimple
Lie algebra $\gs$ is bounded from below by the rank of $\gs$. The same is true for
a nilpotent element. In fact, by a theorem of  de Siebenthal, Dynkin and Kostant \cite[Thm.\ 4.1.6]{CM}, 
the dimension of the centralizer of a nilpotent element 
in a semisimple Lie algebra $\gs$ is bounded from below by the rank of $\gs$ \cite{CM}. 
This proves that $\mathrm{rk}\, \gs =1$ and $\gg = \mathfrak{u}(2)$ or $\gg = \mathfrak{gl}(2,\bR )$, 
since $\dim \gz \le 2$ and $\dim \gg$ is even. 
 \qed

\end{document}